\documentclass[a4paper,10pt,leqno]{amsart}
\date{March 27, 2020}
\usepackage{graphicx} 
\usepackage{verbatim,enumerate}
\usepackage[dvipsnames]{xcolor}
\title[Hypersurfaces with Light-Like Points]{%
  Hypersurfaces with Light-Like Points \\
  in a Lorentzian Manifold II
}
\author[M. Umehara and K. Yamada]{
        M.~Umehara and   
        K.~Yamada        
}   

\address[Umehara]{%
   Department of Mathematical and Computing Sciences,
   Tokyo Institute of Technology,
   Tokyo 152-8552, Japan
}
\email{umehara@is.titech.ac.jp}

\address[Yamada]{%
   Department of Mathematics,
   Tokyo Institute of Technology,
   Tokyo 152-8551, Japan
}
\email{kotaro@math.titech.ac.jp}


\thanks{
The second author
was partially supported by 
the Grant-in-Aid for Scientific Research (B) No.\  
17H02839 from Japan Society for the 
Promotion of Science.
}

\usepackage{amsthm}
\theoremstyle{plain}
 \newtheorem{introtheorem}{Theorem}

 \newtheorem{introfact}[introtheorem]{Fact}
 \newtheorem{theorem}{Theorem}

\theoremstyle{definition}
 
\theoremstyle{remark}
 \newtheorem{remark}[theorem]{Remark}
 \newtheorem*{remark*}{Remark}
 
\newtheorem*{acknowledgements}{Acknowledgements}

\numberwithin{equation}{section}

\newcommand{\mb}[1]{\vect{#1}}
\newcommand{\mc}[1]{\mathcal{#1}}

\newcommand{\vect}[1]{\boldsymbol{#1}}
\newcommand{\R}{\boldsymbol{R}}

\newcommand{\XX}{\mathcal{X}}

\renewcommand{\phi}{\varphi}

\begin{document}
\maketitle

\begin{abstract}
 In the authors' previous work, 
 it was shown that if a
 zero mean curvature  $C^4$-differentiable hypersurface 
 in an arbitrarily given Lorentzian manifold 
 admits a degenerate light-like point,
 then the hypersurface contains a light-like geodesic
 segment passing through the point.
 The purpose of  this paper is to point out that 
 the same conclusion holds with just
 $C^3$-differentiability of the hypersurfaces.
\end{abstract}

\section*{Introduction} \label{sec:1} 
This paper is a remark on authors' previous work
\cite{UY}.
We denote by $M$ a $C^\infty$-differentiable Lorentzian $(n+1)$-manifold
($n\ge 2$). 
Let $U$ be a domain of $(\R^n;u_1,\dots,u_n)$ and
$o(\in U)$ an arbitrary fixed point.
We let $F:U\to M$ be a $C^r$-immersion ($r\ge 3$).
A point $q\in U$ which is neither time-like nor
space-like is called a \emph{light-like point},
which is a zero of the function $B_F$ defined in \cite[(1.2)]{UY}.
A light-like point $q$ is called \emph{degenerate}
if the exterior derivative of $B_F$ vanishes at $q$.
In \cite{UY},
a class $\XX^{r}_0(M,\hat o)$ of germs of 
$C^r$-immersion at $\hat o:=F(o)$ is defined. 
Zero mean curvature immersions in $M$
locally belong to this class. It was shown that:

\begin{introfact}[{\cite[Theorem A]{UY}}] \label{fact:A}
 Suppose that $F\in \XX^{r}_0(M,\hat o)$ 
 is a germ of a $C^r$-immersion $(r\ge 4)$
 at a degenerate light-like point $o$.
 Then the image of $F$ contains a light-like geodesic 
 segment in $M$ 
 passing through $\hat o(=F(o))$
 consisting of only degenerate light-like points.
\end{introfact}

When $M=\R^3_1$ and $F$ is a zero mean curvature
immersion, Fact~\ref{fact:A} was shown  by Klyachin \cite{Kl} 
under the assumption that $F$ is $C^3$-differentiable.
So it is expected that the fact also holds for any 
$C^3$-differentiable $F$.
The above fact was applied to
prove a Bernstein-type theorem in \cite{AHUY}
for entire zero mean curvature graphs without time-like points.
So by the following theorem, the main theorem  
of \cite{AHUY} holds under $C^3$-differentiability  
of the graphs: 

\begin{introtheorem}\label{thm:main}
 In Fact~\ref{fact:A}, the same conclusion holds
 even if $F$ is a germ of a $C^3$-immersion.
 In particular, if $F$ is a $C^3$-differentiable 
 zero mean curvature hypersurface in $M$
 containing a degenerate light-like point $o$, 
 then $F$ contains a light-like geodesic 
 segment passing through $F(o)$
 consisting of only degenerate light-like points.
\end{introtheorem}

In \cite{UY}, we proved Fact~\ref{fact:A} by showing 
the local Lipschitz property of a certain system of ordinary 
differential equations (cf.\ \cite[(4.8)]{UY}).
However, by a careful tracing of the proof of 
\cite[Theorem A]{UY}, we obtain the same conclusion 
under merely $C^3$-differentiability of $F$, as in
the discussions in the next section.
It should be remarked that 
the authors do not know whether the theorem holds under 
$C^2$-differentiability of $F$ or not, which should remain 
an open problem.

\section{Proof of the theorem}
Since $o\in U$ is a light-like point, there exists a
non-zero tangent vector $\vect{v}$ of $U$ at $o$
such that $dF(\vect{v})$ is a light-like vector.
We let $\sigma$ be a light-like geodesic passing through $F(o)$
such that $dF(\vect{v})$ points in the 
tangential direction of $\sigma$ at $F(o)$.
We then take a Fermi-coordinate system 
$(x_0,\dots,x_n)$ of signature $(-+\cdots+)$ centered at 
$F(o)\in M$ along $\sigma$ 
 (cf.\ \cite[Appendix A]{UY}).
When $M$ is the Lorentz-Minkowski space
$\R^{n+1}_1$,  the Fermi coordinate system $(x_0,\dots,x_n)$ can be
taken as the canonical coordinate system.

Using this coordinate system, $\sigma$ can be expressed as
\[
    \sigma(t)=(t,0,\dots,0,t).
\]
Moreover, all of the Christoffel symbols of the 
Lorentzian metric $g$ of $M$
vanish along $\sigma$.
We can express $F$ as (cf.\ \cite[(4.1)]{UY})
\begin{align*}
F(u_1,\dots,u_n)&=\Big(f(x_1(u_1,\dots,u_n),\dots,x_n(u_1,\dots,u_n)),\\
& \hspace{0.2\linewidth}x_1(u_1,\dots,u_n),\dots,x_n(u_1,\dots,u_n)\Big),
\end{align*}
where $f(x_1,\dots,x_n)$ is a $C^3$-function of 
variables $x_1,\dots,x_n$ defined on a neighborhood of the origin $o\in
\R^n$ ($o$ corresponds to the degenerate light-like point of $F$).
In this expression of $F$, $(x_1,\dots,x_n)$ can be considered 
as a local coordinate system of $U$ at $o$ without loss of generality.
We set $y:=x_n$ and use the notation ${}'=d/dy$, 
and set
\begin{equation}\label{eq:a-b}
   a(y):=f(0,\dots,0,y),\quad\text{and}\quad
   b_i(y):=f_{x_i}(0,\dots,0,y)\qquad (i=1,...,n-1),
\end{equation}
where $f_{x_i}:=\partial f/\partial x_i$.
Then $a$ (resp.\ $b_i$ for each $i=1,\dots,n-1$) is a $C^3$
(resp.\ $C^2$)-differentiable function of single variable $y$.
Using these, we define a function $c$ by (cf.\ \cite[(4.2)]{UY}) by
\begin{equation}\label{eq:f}
 f(x_1,\dots,x_{n-1},y)
   =a(y)+\left(\sum_{i=1}^{n-1}b_i(y)x_i\right)+c(x_1,\dots,x_{n-1},y).
\end{equation}
Then $c$ is a function of class $C^2$, and the third derivatives
\[
    c_{yx_ix_j}=c_{x_iyx_j}=f_{yx_ix_j}, \qquad c_{yyx_i}=f_{yyx_i}-b_i'
\]
exist for $i,j=1,\dots,n-1$, which are continuous functions.
By \eqref{eq:a-b} and \eqref{eq:f}, one can easily show that
\begin{align}\label{eq:N1}
 c(0,\dots,0,y)&=c_{x_i}(0,\dots,0,y)
    = c_{y}(0,\dots,0,y)\\
 &=c_{yx_i}(0,\dots,0,y)=
       c_{yyx_i}(0,\dots,0,y)=0
\nonumber
\end{align}
and that
\begin{equation}\label{eq:N2}
\begin{aligned}
   f_{x_ix_jx_k}(0,\dots,0,y)&=
   c_{x_i x_j x_k}(0,\dots,0,y),\\
   f_{yx_ix_j}(0,\dots,0,y)&=
   c_{yx_ix_j}(0,\dots,0,y)=
   c_{x_iyx_j}(0,\dots,0,y),\\
   f_{yyx_i}(0,\dots,0,y)&= b''_{i}(y)
\end{aligned}
\end{equation}
for $i,j,k=1,\dots,n-1$.

\begin{remark}\label{rem:X}
 In \cite{UY} and \cite{AHUY}, we used
 the function $h_{ij}$ satisfying
 \begin{equation}\label{eq:f2}
  f(x_1,\dots,x_{n-1},y)
   =a(y)+\sum_{i=1}^{n-1}b_i(y)x_i+\sum_{i,j=1}^{n-1} 
   \frac{h_{ij}(x_1,\dots,x_{n-1},y)}2 x_ix_j.
 \end{equation}
 Moreover, we wrote that
 $h_{ij}$ are $C^1$-functions. 
 However it may not be true, in general.
 We should avoid this expression. 
 In fact, for example, even if
 $h_{ii}$ is a $C^1$-function, 
 the following computation is not allowed
 \[
 \Big(h_{ii} x_i^2\Big)_{x_ix_i}
 =(h_{ii})_{x_ix_i}x_i^2+2(h_{ii})_{x_i}x_i+2h_{ii}.
 \]
 So it should be replaced
 \eqref{eq:f2} by  \eqref{eq:f} 
 in \cite{UY} and the final remark of
 \cite{AHUY}. 
\end{remark}

The condition that $o$ is a degenerate light-like point
is written as (cf.\ \cite[(4.3) and (4.5)]{UY})
\begin{equation}\label{eq:Ini}
 a(0)=0,\quad a'(0)=1,\quad b_i(0)=b'_i(0)=0 \qquad (i=1,\dots,n-1),
\end{equation}
where ${}'=d/dy$.

Since $F\in \XX^{3}_0(M,\hat o)$, there
exists a $C^1$-function $\phi$ such that  
(cf.\ \cite[(3.1) and (3.3)]{UY})
\begin{equation}\label{eq:tildeA}
  \tilde A:=A-\phi B
\end{equation}
vanishes identically on $U$,
where 
\[
    B:=B_F,\qquad A:=A_F
\]
are functions associated with $F$
defined in \cite[(2.4) and (2.6)]{UY}.
We remark that $B,\nabla B$ and $A$ can be
expressed without use of derivatives of the Christoffel symbols,
where $\nabla B:=(B_{x_1},\dots,B_{x_n})$.
Since $(x_0,\dots,x_n)$ is a Fermi coordinate system,
the expressions of $B$, $\nabla B$ and $A$
along $\sigma$ coincide with
those of $\R^{n+1}_1$ given in \cite[Appendix B]{UY}.
For the Lorentzian metric $g$ of $M$, set
\[
  g_{\alpha\beta}^{}
               :=g(\partial_{x_\alpha},\partial_{x_\beta}),\quad 
 \hat g_{\alpha\beta}^{}:=g_{\alpha\beta}^{}
          \circ F \qquad (\alpha,\beta=0,\dots,n),
\]
where $\partial_{x_\alpha}:=\partial/\partial x_\alpha$.
Since we can write
\begin{equation}\label{eq:Fi}
  F_{x_i}=f_{x_i}\partial_0+\partial_{x_i} \qquad (i=1,\dots,n),
\end{equation}
the coefficients $s_{i,j}^{}$ of the induced metric $ds^2$
with respect to the local coordinate system $(U;x_1,\dots,x_n)$ 
given in  \cite[(2.3)]{UY} can be computed as
\[
   s_{i,j}^{}=f_{x_i}f_{x_j}\hat g_{00}^{}+
   f_{x_i}\hat g_{0j}^{}+f_{x_j}\hat g_{i0}^{}+
        \hat g_{ij}^{} \qquad (i,j=1,\dots,n).
\]
Since $\{\hat g_{\alpha\beta}^{}\}_{\alpha,\beta=0}^n$ are functions of $f$,
each $s_{i,j}^{}$ is also a function of $f$ and $f_I$, 
where
\[
   f_I:=\{f_{x_i}\}_{i=1}^n.
\]
In particular,  $B$ is a function of $f$ and $f_I$.
Also each component $\tilde s^{i,j}$ ($i,j,=1,\dots,n$)
of the cofactor matrix of the $n\times n$-matrix  $(s_{i,j}^{})$
is also a function of $f$ and $f_I$.

By \eqref{eq:Fi}, $\mu_{i}^\beta$ ($i=1,\dots,n$, $\beta=0,\dots,n$)
in \cite[Page 3410]{UY} are given by
\[
  \mu_{i}^\beta=
      \begin{cases}
       1 & \text{if $\beta=i$}, \\
       f_{x_i} & \text{if $\beta=0$}, \\
       0 &  \text{if $\beta\ne 0,i$}, 
      \end{cases}
\]
and
\[
   \mu_{i,\beta}^{}=\sum_{\alpha=0}^n \mu_{i}^\alpha \hat g_{\alpha\beta}
   \qquad
   (i=1,\dots,n,\,\,\beta=0,\dots,n)
\]
are functions of $f$ and $f_I$.
In particular, each component $\tilde \nu_i$ 
of the normal vector field $\tilde \nu$ (cf.\ \cite[(2.5)]{UY})
is also a function of $f$ and $f_I$.
We denote by $\{\Gamma_{\alpha\beta}^\gamma\}_{\alpha,\beta,\gamma=0}^n$
the Christoffel symbols of $g$ with respect to 
the coordinate system $(x_0,\dots,x_n)$ of $M$,
and set
\[
  \hat \Gamma_{\alpha\beta}^\gamma:=\Gamma_{\alpha\beta}^\gamma
  \circ F
  \qquad (\alpha,\beta,\gamma=0,\dots,n),
\]
each of which can be considered as a function of $f$.
Since $F_{x_j}=f_{x_j}\partial_{x_0}+\partial_{x_j}$,
we can write
\[
  D_{\partial_{x_i}}F_{x_j}
    =f_{x_ix_j}\partial_{x_0}+
    \sum_{\alpha=0}^n 
    \left(
         f_{x_i}f_{x_j}\hat \Gamma_{00}^\alpha+
	 f_{x_i}\hat \Gamma_{0j}^\alpha+
	 f_{x_j}\hat \Gamma_{i0}^\alpha+
	 \hat \Gamma_{ij}^\alpha
	 \right)\partial_{x_\alpha},
\]
where $D$ is the Levi-Civita connection induced by $g$.
In particular, by \cite[(2.6)]{UY}, $A$ is also a function of 
$f$, $f_I$, $f_{IJ}$
and $f_{In}$, where
\[
  f_{IJ}:=\{f_{x_ix_j}\}_{i,j=1}^{n-1},\quad
  f_{In}:=\{f_{x_ix_n}\}_{i=1}^{n-1}.
\]
As a consequence, $\tilde A$ (cf.\ \eqref{eq:tildeA})
has the following expression:
\begin{equation}\label{eq:A}
 \tilde A=\mc R(f, f_I)f_{x_nx_n}+\mc S(f_I,f_{IJ},f_{In},\phi),
\end{equation}
where $\mc R$ is a function of $f$, $f_I$, 
and $\mc S$
is a function of $f$, $f_{I}$, $f_{IJ}$, $f_{In}$ and $\phi$.
For example, if $M=\R^{n+1}_1$ and $\phi=0$, then $\tilde A=A$ and
\begin{align}
\mc R&=1-\sum_{j=1}^{n-1}f_{x_j}^2, \\
\mc S&=\sum_{k=1}^{n-1}
         \left(1-f_{x_n}^2-\sum_{j=1, j\neq k}^{n-1}
                     f_{x_j}^2\right)f_{x_kx_k}\\ \nonumber
    &\hphantom{==+}    + 2\sum_{1\leq j<k\leq n-1} f_{x_j} f_{x_k} 
            f_{x_jx_k}
        + 2\sum_{j=1}^{n-1} f_{x_j} f_{x_n} f_{x_jx_n}.
\end{align}
Although $f_{I}$ (resp.\ $f_{IJ}$ and $f_{In}$)  
has (resp.\ have) only $C^2$-differentiability (resp.\ $C^1$-differentiability),
$\mc R$ and $\mc S$ depend smoothly on
variables $f,f_{I},f_{IJ},f_{In}$ and $\phi$, because
$M$ itself is a $C^\infty$-manifold.
Differentiating \eqref{eq:A},
we have the following expression
of the derivative $\tilde A_{x_i}$ with respect to $\partial_{x_i}$:
\[
  \tilde A_{x_i}=\mc R(f, f_I)f_{x_nx_nx_i}+\mc T(f_I,f_{IJ},f_{In},
  f_{IJK},f_{IJn},\phi, \phi_I),
\]
where $\mc T$
is a function of $f_I$,
$f_{IJ}$, $f_{In}$, $f_{IJK}$, $f_{IJn}$, $\phi$, $\phi_I$ and
\[
   f_{IJK}:=\{f_{x_ix_jx_k}\}_{i,j=1}^{n-1},\quad
   f_{IJn}:=\{f_{x_nx_ix_j}\}_{i,j=1}^{n-1}
\]
For example, if $n=2$, 
$\mc R$ is a function of $f,f_{x_1}$ 
and $\mc S$ is a function of $\phi$, $f$, $f_{x_1}$, 
$f_{x_1x_1}$, $f_{x_1x_2}$.
So we have the following:
\begin{align*}
&A_{x_1}=\mc R_{f}(f, f_{x_1}) f_{x_1} f_{x_2x_2}+\mc R_{f_{x_1}}
(f, f_{x_1}) f_{x_1x_1} f_{x_2x_2} 
+\mc S_{\phi}(f,f_{x_1},f_{x_1x_1},f_{x_1x_2},\phi)\phi_{x_1}	\\
&\phantom{aa}+
\mc S_f(f,f_{x_1},f_{x_1x_1},f_{x_1x_2},\phi) f_{x_1}+
\mc S_{f_{x_1}}(f,f_{x_1},f_{x_1x_1},f_{x_1x_2},\phi) f_{x_1x_1}\\
&\phantom{aaa}+\mc S_{f_{x_1x_1}}(f,f_{x_1},
f_{x_1x_1},f_{x_1x_2},\phi) f_{x_1x_1x_1}+
\mc S_{f_{x_1x_2}}(f,f_{x_1},f_{x_1x_1},f_{x_1x_2},\phi) 
f_{x_1x_2x_1},
\end{align*}
where, for example,
\[
  \mc R_{f}(f, f_{x_1}):=
   \frac{\partial \mc R(f, f_{x_1})}{\partial f},\qquad
  \mc R_{f_{x_1}}(f, f_{x_1}):=\frac{\partial \mc R(f, f_{x_1})}
   {\partial f_{x_1}}.
\]

Since $\tilde A$ vanishes identically, we have
$\tilde A=0$ and $\tilde A_{x_i}=0$ for
$i=1,\dots,n-1$.
We set $\mb x:=(x_1,\dots,x_{n-1})$
and substitute $\mb x=(0,\dots,0)$
to them.
Then they are functions of one variable $y$,
and we have
\begin{align}\label{eq:E1}
 &(\mathfrak a:=)\tilde A|_{\mb x=(0,\dots,0)}=0,\\
 &(\mathfrak a_i:=)\tilde A_{x_i}|_{\mb x=(0,\dots,0)}=0
 \qquad (i=1,\dots,n-1).
\label{eq:E2}
\end{align}
Although $\mathfrak a$ and $\mathfrak a_i$ ($i=1,\dots,n-1$) are
all functions of $y$, we now would like to
think that they are functions of
$a$, $a'$, $a''$, $b_I$, $b'_I$, $b''_I$, $c_{IJ}$, $c_{IJK}$ and $\phi$,
where 
\begin{align*}
 & c_{ij}(y):=c_{x_ix_j}(0,\dots,0,y),\quad
 c_{ijk}(y):=c_{x_ix_jx_k}(0,\dots,0,y),\\
 &\hat\phi(y):=\phi(0,\dots,0,y),\quad
 \hat\phi_i^{}(y):=\phi_{x_i}(0,\dots,0,y)
\end{align*}
and
\begin{align*}
 &b_I:=\{b_i\}_{i=1}^{n-1},\quad
 b'_I:=\{b'_i\}_{i=1}^{n-1},\quad
 b''_I:=\{b''_i\}_{i=1}^{n-1},\\
 &\hat\phi_I^{}:=\{\hat \phi_i^{}\}_{i=1}^{n-1},\quad
 c_{IJ}:=\{c_{ij}\}_{i,j=1}^{n-1},\quad 
 c_{IJK}=\{c_{ijk}\}_{i,j,k=1}^{n-1}.
\end{align*}
Then we have
the following expressions 
\begin{align}\label{eq:U1}
 \mathfrak a&=
 \Lambda^1(a,a',b_I,b'_I,c_{IJ},c_{IJK},\hat \phi)a''
 +\Lambda^2(a,a',b_I,b'_I,c_{IJ},c_{IJK},\hat \phi), \\
 \label{eq:U2}
 \mathfrak a_i&=
 \Lambda^3_i(a,a',b_I,b'_I,c_{IJ},c_{IJK},\hat \phi,\hat \phi_I)a''\\
 &\phantom{aa}
 +\Lambda^4_i(a,a',b_I,b'_I,c_{IJ},c_{IJK},\hat \phi,\hat \phi_I)b''_i 
 +\Lambda^5_i(a,a',b_I,b'_I,c_{IJ},c_{IJK},\hat \phi,\hat \phi_I),
 \nonumber
\end{align}
where $\Lambda^1,\Lambda^2$ and $\Lambda^l_i$
($l=3,\dots,5$, $i=1,\dots,n-1$) are functions which 
are $C^\infty$-differentiable with respect to the parameters
$a$, $a'$, $b_I$, $b'_I$, $c_{IJ}$, $c_{IJK}$, $\hat \phi$, $\hat \phi_I$.
For example, the explicit expressions of
$\mathfrak a$ and $\mathfrak a_i$ for $M=\R^{3}_1$
are given in \cite[(4.4) and (4.5)]{Geloma}.

Since $\mc R=1$ at $o$,  
$\Lambda^1$ and $\Lambda^4_i$ ($i=1,\dots,n-1$)
do not vanish at $o$ 
under the initial condition \eqref{eq:Ini}.
Regarding
$c_{IJ}$, $c_{IJK}$, $\hat\phi$ and $\hat\phi_{I}$
are functions of $y$,
we can rewrite 
\eqref{eq:E1} and \eqref{eq:E2}
in the following forms
\begin{align}\label{eq-1a}
a'' &=P(y, a,a',b_I,b'_I), \\
\label{eq-1b}
b''_i &=Q_i(y,a,a',b_I,b'_I)
\qquad (i=1,2,\dots,n-1),
\end{align}
where  $P,Q_i$ 
are functions of $y,a,a',b_I,b'_I$.

For example, we consider the case $M=\R^{n+1}_1$.
If we set
\begin{gather*}
  c_{ij}(y):=c_{x_ix_j}(0,\dots,0,y),\quad
  c_{ijk}(y):=c_{x_ix_jx_k}(0,\dots,0,y),\\
  c'_{ij}(y):=c_{x_ix_jy}(0,\dots,0,y),
\end{gather*}
then
\begin{align}\label{eq:a1}
\alpha&=C a''
      + \sum_{j=1}^{n-1}c_{jj}(D+b_j^2)
      + 2 \sum_{i<j} b_{i}b_{j} c_{ij}
                    + 2 \sum_{i=1}^{n-1}a'b_{i}b'_{i},\\
\label{eq:a2}
\alpha_l
  &=C b_{l}''
    - 2\left( \sum_{j=1}^{n-1}b_{j}c_{jl}\right)a''
    + \sum_{j=1}^{n-1}(D+b_j^2)c_{jjl}\\
    &- 2\sum_{j=1}^{n-1}
      \left(a' b_{l}' 
   + \sum_{i=1,i\neq j}^{n-1}b_{i}c_{il}\right)
      c_{jj}
        + 2 \sum_{i<j}
       \left(b_{j}c_{il}c_{ij}+b_{i}c_{jl}c_{ij}+
             b_{i}b_{j}c_{ijl}\right)
\nonumber
\\
    &+
     2 \sum_{i=1}^{n-1}
       \left(c_{il}a' b'_{i} + b_{i}b_{l}'b_{i}'
                 + b_{i}a'c'_{il}\right)
\nonumber
\end{align}
hold for $l=1,\dots,n-1$, 
where $C:=1-\sum_{i=1}^{n-1} b_{i}^2$ 
and $D:=1-(a')^2-\sum_{i=1}^{n-1}b_{i}^2$.

In particular,
\eqref{eq-1a} and
\eqref{eq-1b}
give a normal form for a system of ordinary differential
equations thinking $a,b_i$ ($i=1,\dots,n-1$)
are unknown functions, 
and can be considered as precise expressions of 
\cite[(4.12)]{UY}.
Since $P,Q_i$ ($i=1,\dots,n-1$)
are $C^\infty$-differentiable with respect to
the variables
$a$,$a'$, $b_I$, $b'_I$,
this
system of ordinary differential equations
satisfies the local Lipschitz condition. 
So the uniqueness of the solution follows.
On the other hand, 
by applying
\cite[Proposition 4.2]{UY},
\begin{equation}\label{eq:C}
a(y)=y,\qquad b_i(y)=0 \qquad (i=1,...,n-1)
\end{equation}
gives a solution of this system satisfying
the initial condition \eqref{eq:Ini}.
However, we should remark that 
in the proof of \cite[Proposition 4.2]{UY},
we set (cf.\ (4.14))
\[
  f_0(x_1,\dots,x_n):=x_n+\sum_{1\le j\le k \le n-1}
           c_{j,k}(x_1,\dots,x_{n})x_jx_k.
\]
However, by the same reason as Remark \ref{rem:X},
we should replace it by
\[
  f_0(x_1,\dots,x_n):=x_n+c(x_1,\dots,x_{n}),
\]
where $c$ is a $C^3$-function satisfying
\eqref{eq:N1} and \eqref{eq:N2}.
Then 
\[
    x_n\mapsto (f_0(0,...,0,x_n),0,...,0,x_n)
\]
gives a light-like geodesic of $M$ consisting of degenerate light-like
points.

So, we can conclude that \eqref{eq:Ini} 
implies \eqref{eq:C}, that is,
we have proven the theorem.

\begin{remark*}
We point out that there is a minor
 typographical error in the proof of Theorem~D 
 in Section 6 in \cite{UY}.
 In fact, in the beginning of the proof, 
 we had set $f(x,y)=a(y)+b(y)x+ h(x,y)x^3$, but it
 should be $f(x,y)=a(y)+b(y)x+ h(x,y)x^2$. The remaining
 arguments can be read without any need for corrections.
\end{remark*}

\begin{acknowledgements}
 The authors thank Professor 
 Udo Hertrich-Jeromin for valuable comments.
\end{acknowledgements}


\end{document}